\documentclass[letterpaper,10pt,conference]{ieeeconf} 

\IEEEoverridecommandlockouts
\let\labelindent\relax
\usepackage{pgfplotstable}
\usepackage{pgfplots}
\usepackage{times,amsmath,epsfig,amssymb,graphicx,amsfonts}
\usepackage{empheq}
\usepackage{graphicx}
\usepackage{psfrag}
\usepackage{fge}
\usepackage{amsmath}
\usepackage{setspace}
\usepackage{color}
\usepackage{amsfonts}   
\usepackage{amssymb}    
\usepackage{mathrsfs}
\usepackage{setspace}
\usepackage{dblfloatfix}
\usepackage{tikz}
\usepackage{tkz-tab}
\usepackage{algorithm}
\usepackage{algcompatible}
\usepackage{lipsum}
\usepackage{color}
\usepackage{mathrsfs}
\usepackage{epstopdf}
\usepackage{enumerate}
\usepackage{enumitem}
\usepackage{url}
\usepackage{cite}
\usepackage{tabularx}
\usepackage[pass]{geometry}
\usepackage{chngcntr}

\floatname{algorithm}{Algorithm}

\floatname{algorithm}{Algorithm}

\newtheorem{mydef}{Definition}

\newtheorem{mytheorem}{Theorem}
\newtheorem{mylemma}{Lemma}
\newtheorem{myremark}{Remark}

\newtheorem{myexample}{Example}
\newtheorem{myproblem}{Problem}

\newcounter{ale}

\newenvironment{liste}{\begin{itemize}}{\end{itemize}}
\newcommand{\aliste}{\begin{liste} \setcounter{ale}{1}}
\newcommand{\zliste}{\end{liste}}

\newcommand{\calA}{{\cal A}}
\newcommand{\calB}{{\cal B}}

\newcommand{\calP}{{\cal P}}

\newcommand{\calV}{{\cal V}}

\newcommand{\calY}{{\cal Y}}

\title{\LARGE \bf Resilient Monotone Submodular Maximization}
\author{Vasileios Tzoumas,{$^{1}$} Konstantinos Gatsis,{$^{1}$} Ali Jadbabaie,{$^{2}$} George J.~Pappas{$^{1}$}
\thanks{$^{1}$The authors are with the Department of Electrical and Systems Engineering, University of Pennsylvania, Philadelphia, PA 19104-6228 USA (email: {\fontsize{8}{8}\selectfont\ttfamily\upshape \{vtzoumas, kgatsis, pappasg\}@seas.upenn.edu}).}
\thanks{$^{2}$The author is with the Institute for Data, Systems and Society, Massachusetts Institute of Technology, Cambridge, MA 02139 USA (email: {\fontsize{8}{8}\selectfont\ttfamily\upshape jadbabai@mit.edu}).}
\thanks{This work was supported in part by TerraSwarm, one of six centers of STARnet, a Semiconductor Research Corporation program sponsored by MARCO and DARPA, and in part by AFOSR Complex Networks Program.}
}

\begin{document}
\maketitle

\begin{abstract}
In this paper, we focus on applications in machine learning, optimization, and control that call for the resilient selection of a few elements, e.g. features, sensors, or leaders, against a number of adversarial denial-of-service attacks or failures. 
\mbox{In general}, such resilient optimization problems are hard and cannot be solved exactly in polynomial time, even though they may involve objective functions that are monotone and submodular. 
In this paper, we provide for the solution of such optimization problems the first scalable approximation algorithm that is valid for any number of attacks or failures and which, for functions with low curvature, guarantees superior approximation performance. Notably, the curvature has been known to tighten approximations for several non-resilient optimization problems, yet its effect on resilient optimization had hitherto been unknown. 
We complement our theoretical analyses with empirical evaluations.
\end{abstract}

\section{Introduction}\label{sec:Intro}

During the last decade, researchers  in machine learning, optimization, and control have focused on questions such as:
\begin{itemize}[leftmargin=*]
\item \textit{(Sensor selection)} How many sensors do we need to deploy in a large water distribution network to detect a contamination outbreak as fast as possible?~\cite{krause08efficient}
\item \textit{(Feature selection)} Which few features do we need to select from the data flood to optimally train a spam e-mail classifier?~\cite{das2011spectral}
\item \textit{(Leader selection)} Which UAVs in a multi-UAV system do we need to choose as leaders so the system can complete a surveillance task despite communication noise?~\cite{clark2016submodularity}
\end{itemize}
The effort to answer such questions has culminated in a plethora of papers on topics such as actuator placement for controllability~\cite{olshevsky2014minimal, bullo2014,sergio2015minimal,summers2014submodularity, tzoumas2016minimal,nozari2016time}; 
sensor scheduling for target tracking~\cite{krause2008near,jawaid2015submodularity,sharma2015greedy,tzoumas2016near,tzoumas2016scheduling,zhang2017kalman}; and visual cue selection for robotic navigation~\cite{carlone2016attention,iyer2013fast}.  
Notably, in all aforementioned papers the underlying optimization problem is of the form
\begin{equation}\label{eq:non_res}
\max_{\mathcal{A}\subseteq \mathcal{V}, |\mathcal{A}|\leq \alpha} \; \; f(\mathcal{A}),
\end{equation}
where the set function $f$ exhibits monotonicity and submodularity, a diminishing returns property~\cite{krause08efficient,das2011spectral,clark2016submodularity,olshevsky2014minimal, bullo2014,sergio2015minimal,sharma2015greedy,summers2014submodularity, tzoumas2016minimal,nozari2016time,jawaid2015submodularity,tzoumas2016near,tzoumas2016scheduling,zhang2017kalman,iyer2013fast,krause2008near,carlone2016attention}.
\mbox{In words}, Problem~\eqref{eq:non_res} aims to find a set $\mathcal{A}$ of $\alpha$ elements from the finite ground set $\mathcal{V}$, such that $\mathcal{A}$ maximizes~$f$.  This problem is NP-hard, yet several good approximation algorithms have been proposed for its solution, such as the greedy~\cite{nemhauser1978best}.

But sensors and actuators fail~\cite{krause2008robust}; features can become obsolete~\cite{globerson2006nightmare}; and leaders can be attacked~\cite{clark2012leader}. \mbox{For example},
sensors may fail due to malfunctions or adversarial attacks.
In~such~scenarios, questions such as the following arise:
\begin{itemize}[leftmargin=*]
\item Where to place a few actuators in a system, when some of them may fail?~\cite{pequitoJ7}
\item Which sensors to activate to track an adversarial target that can jam a fraction of the activated sensors?~\cite{das2008sensor,laszka2015resilient}
\item Or, which features to select to train a robust machine learning model to changing features?~\cite{xu2006jamming}
\end{itemize} 
In such scenarios, the optimization problem we need to address takes the form
\begin{equation}\label{eq:res}
\max_{\mathcal{A}\subseteq \mathcal{V}, |\mathcal{A}|\leq \alpha} \; \min_{\mathcal{B}\subseteq \mathcal{A}, |\mathcal{B}|\leq \beta}\; \; f(\mathcal{A}\setminus\mathcal{B}),
\end{equation}
where $\beta\leq \alpha$, which is a resilient formulation of Problem~\eqref{eq:non_res}. In words, Problem~\eqref{eq:res} aims to find a set $\mathcal{A}$ of $\alpha$ elements such that $\mathcal{A}$ is resilient against the worst possible removal of $\beta$ of its elements.
Importantly, this formulation is suitable when we have no prior on the failure or attack mechanism.  

The most relevant papers on the resilient monotone submodular optimization Problem~\eqref{eq:res} are~\cite{krause2008robust} and~\cite{Orlin2016}. \mbox{In particular}, Problem~\eqref{eq:res} was introduced in~\cite{krause2008robust}, where the authors proposed an approximation algorithm for the more general problem $\max_{\mathcal{A}\subseteq \mathcal{V}, |\mathcal{A}|\leq \alpha} \min_{i\in\{1,2,\ldots,m\}}f_i(\mathcal{A})$, where each $f_i$ is monotone submodular. In more detail, the algorithm in~\cite{krause2008robust} guarantees a high value for Problem~\eqref{eq:res}
by allowing sets $\mathcal{A}$ up to size $\alpha [1+2\log(\alpha\beta\log(|\mathcal{V}|))]$, instead of $\alpha$. Nonetheless, it runs with $O(|\mathcal{V}|^2(\frac{\alpha}{\beta})^\beta$) evaluations of $f$, which is exponential in the number of possible removals $\beta$, and quadratic in the number of available elements for resiliency $|\mathcal{V}|$. This limits its applicability in large-scale settings where both $\beta$ and $|\mathcal{V}|$ can be in the order of several thousands~\cite{shoukry2016smt}. In~\cite{Orlin2016}, and its arxiv version~\cite{orlin2015robust}, the authors prove that Problem~\eqref{eq:res} is NP-hard. In addition, they provide an algorithm for Problem~\eqref{eq:res} that runs with $O(|\mathcal{V}|(\alpha-\beta))$ evaluations of $f$, and which, for $\alpha, \beta \longrightarrow +\infty$, guarantees an approximate value at least $\simeq 29\%$ the optimal, the first approximation performance bound for Problem~\eqref{eq:res}. However, in~\cite{Orlin2016} the proposed algorithm is valid only when the number $\beta$ of possible failures and attacks is up to $\sqrt{2\alpha}$, whereas in practice the number $\beta$ of failures and attacks can be of the same order as the number $\alpha$ of selected sensors, actuators, etc. For example, it may be the case that $\beta$ is up to $\alpha/2$~\cite{chong2015}.

In this paper, we show how the notion of the curvature ---deviation from modularity (additivity)--- of a function can be used to provide a scalable algorithm for the resilient maximization Problem~\eqref{eq:res} that has maximum resiliency, and for submodular functions with low curvature, superior approximation performance. 
Notably, low curvature submodular functions are involved in a series of applications~\cite{jegelka2012combinatorial,sharma2015greedy}, such as sensor placement for mutual information maximization~\cite{krause2008near}; feature selection for Gaussian process regression~\cite{sharma2015greedy}; and active learning for speech recognition~\cite{iyer2013fast}.

In more detail, exploiting the curvature of the monotone submodular function $f$, denoted henceforth by $\kappa_f$, in the resilient maximization Problem~\eqref{eq:res}, we~provide an algorithm (Algorithm~\ref{alg:rob_sub_max})
with the properties:
\begin{itemize}[leftmargin=*]
\item Algorithm~\ref{alg:rob_sub_max} is valid for any number of selections for resiliency $\alpha$ and any number of failures or attacks $\beta$;
\item Algorithm~\ref{alg:rob_sub_max} runs with $O(|\mathcal{V}|(\alpha-\beta))$ evaluations of $f$;
\item Algorithm~\ref{alg:rob_sub_max} guarantees the approximation performance bound $$\frac{f_\text{Algorithm~1}}{f^\star}\geq \max\left(1-\kappa_f, \frac{1}{\beta+1}\right)\frac{1}{\kappa_f}\left(1-e^{-\kappa_f}\right),$$
where $f^\star$ is the (optimal) value of Problem~\eqref{eq:res}, and $f_\text{Algorithm~1}$ is the (approximate) value achieved by Algorithm~1 for Problem~\ref{eq:res}.  Notably, the notion of curvature we use for the monotone submodular function $f$ is such that the curvature $\kappa_f$ takes only the values $0\leq \kappa_f\leq 1$.
\end{itemize}

Overall, Algorithm~\ref{alg:rob_sub_max} improves upon the state-of-the-art algorithms for Problem~\eqref{eq:res} as follows:
\begin{itemize}[leftmargin=*]

\item \textit{High resiliency:} Algorithm~\ref{alg:rob_sub_max} is the first scalable algorithm for Problem~\eqref{eq:res} with bounded approximation performance for any number of failures or attacks $\beta \leq \alpha$.

\item \textit{High approximation performance:} For low curvature values ($\kappa_f\leq 0.71$), Algorithm~\ref{alg:rob_sub_max} is the first scalable algorithm to exhibit approximation performance for Problem~\eqref{eq:res} at least $29\%$ the optimal.

\vspace{0.8mm}
For example, for the central problem in machine learning of Gaussian process regression with RBF kernels~\cite{krause2008near,bishop2006pattern}, Algorithm~\ref{alg:rob_sub_max} guarantees almost exact approximation performance (approximate value $\simeq 100\%$ the optimal).
	\end{itemize}



\section{Resilient Submodular Maximization}\label{sec:resilient_sub_max}
 
We state the resilient maximization problem considered in this paper.  To this end, we start with the definitions of monotone and submodular set functions.

\textit{Notation:} 
For any set function $f:2^\mathcal{V}\mapsto \mathbb{R}$ on a ground set $\mathcal{V}$, and any element $x\in \mathcal{V}$, $f(x)$ denotes $f(\{x\})$.
 
\begin{mydef}[Monotonicity]
Consider any finite ground set $\mathcal{V}$. The set function $f:2^\mathcal{V}\mapsto \mathbb{R}$ is non-decreasing if and only if for any $\mathcal{A}\subseteq \mathcal{A}'\subseteq\mathcal{V}$, $f(\mathcal{A})\leq f(\mathcal{A}')$.\hfill $\blacktriangle$
\end{mydef}

In words, a set function $f:2^\mathcal{V}\mapsto \mathbb{R}$ is non-decreasing if and only if adding more elements in any set $\mathcal{A}\subseteq\mathcal{V}$ cannot decrease the value of $f(\mathcal{A})$.
 
\begin{mydef}[Submodularity~{\cite[Proposition 2.1]{nemhauser78analysis}}]\label{def:sub}
Consider any finite ground set $\mathcal{V}$.  The set function $f:2^\mathcal{V}\mapsto \mathbb{R}$ is submodular if and only if
\begin{itemize}
\item for any sets $\mathcal{A}\subseteq \mathcal{V}$ and $\mathcal{A}'\subseteq \mathcal{V}$, it is $f(\mathcal{A})+f(\mathcal{A}')\geq f(\mathcal{A}\cup \mathcal{A}')+f(\mathcal{A}\cap \mathcal{A}')$;
\item equivalently, for any sets  $\mathcal{A}\subseteq \mathcal{A}'\subseteq\mathcal{V}$ and any element $x\in \mathcal{V}$, it is $f(\mathcal{A}\cup \{x\})-f(\mathcal{A})\geq f(\mathcal{A}'\cup \{x\})-f(\mathcal{A}')$.

\hfill $\blacktriangle$
\end{itemize}
\end{mydef}

In words, a set function $f:2^\mathcal{V}\mapsto \mathbb{R}$ is submodular if and only if it satisfies the following intuitive diminishing returns property: for any element $x\in \mathcal{V}$, the marginal gain $f(\mathcal{A}\cup \{x\})-f(\mathcal{A})$ diminishes as the set $\mathcal{A}$ grows; equivalently, for any set $\mathcal{A}\subseteq \mathcal{V}$ and any element $x\in \mathcal{V}$, the marginal gain $f(\mathcal{A}\cup \{x\})-f(\mathcal{A})$ is non-increasing.

In this paper, we consider the problem of resilient monotone submodular maximization, defined as follows.
\begin{myproblem}\label{pr:robust_sub_max} 
Consider
\begin{itemize}
\item a finite ground set $\mathcal{V}$;
\item a submodular and monotone set function $f:2^\mathcal{V} \mapsto \mathbb{R}$ such that (without loss of generality) $f$ is non-negative and $f(\emptyset)=0$;
\item and two integers $\alpha$ and $\beta$ such that $0\leq\beta \leq \alpha\leq |\mathcal{V}|$.
\end{itemize}
The problem of resilient monotone submodular maximization is to select a set $\mathcal{A}$ of $\alpha$ elements from the ground set~$\mathcal{V}$, such that $f(\mathcal{A})$ is resilient against the worst possible removal~$\mathcal{B}$ of $\beta$ of $\mathcal{A}$'s elements. In mathematical notation,
\begin{equation*}
\underset{\mathcal{A} \subseteq \mathcal{V}, |\mathcal{A}| \leq \alpha}{\max} \; \; \underset{\mathcal{B} \subseteq \mathcal{A}, |\mathcal{B}| \leq \beta}{\min} \; \; f(\mathcal{A}\setminus \mathcal{B}).
\end{equation*}

\vspace{-3mm}
\hfill $\blacktriangle$
\end{myproblem}

The resilient maximization Problem~\ref{alg:rob_sub_max} may be interpreted as a two-stage perfect information sequential game~\cite[Chapter~4]{myerson2013game}, where the player that plays first chooses a set~$\mathcal{A}$, and the player that plays second, knowing $\mathcal{A}$, chooses a set-removal~$\mathcal{B}$ from $\mathcal{A}$.

\section{Algorithm for resilient submodular maximization} \label{sec:main}

We present the first scalable algorithm for the resilient maximization Problem~\ref{pr:robust_sub_max}, and show that this algorithm is valid for any number of failures and attacks, and that for functions with low curvature it guarantees superior approximation performance.  We begin by presenting the definition of curvature of monotone submodular functions. 

\subsection{Curvature of monotone submodular functions}\label{subsec:curv}

We define the curvature of monotone submodular functions, which we use to quantify the approximation performance of the proposed algorithm in this paper. To this end, we start with the definition of modular set functions.

\begin{mydef}[Modularity]\label{def:modular}
Consider any finite ground set~$\mathcal{V}$.  The set function $f:2^\mathcal{V}\mapsto \mathbb{R}$ is modular if and only if for any set $\mathcal{A}\subseteq \mathcal{V}$, it is $f(\mathcal{A})=\sum_{v\in \mathcal{A}}f(v)$. \hfill $\blacktriangle$
\end{mydef}
\vspace{0.3mm}

In words, a set function $f:2^\mathcal{V}\mapsto \mathbb{R}$ is modular if through $f$ all elements in the ground set $\mathcal{V}$ cannot substitute each other, since Definition~\ref{def:modular} implies that for any set $\mathcal{A}\subseteq\mathcal{V}$ and any element $x\in \mathcal{V}\setminus\mathcal{A}$, it is $f(\{x\}\cup\mathcal{A})-f(\mathcal{A})= f(x)$; that is, in the presence of $\mathcal{A}$, $x$~retains its contribution to the value of $f(\{x\}\cup\mathcal{A})$.
In contrast, for a submodular set function $g:2^\mathcal{V}\mapsto \mathbb{R}$, the~elements in $\mathcal{V}$ \emph{can} substitute each other, since Definition~\ref{def:sub} implies $g(\{x\}\cup\mathcal{A})-g(\mathcal{A})\leq g(x)$; that is, in the presence of $\mathcal{A}$, $x$ may \mbox{lose part of its contribution to the value of~$g(\{x\}\cup\mathcal{A})$}.

\begin{mydef}[Curvature]\label{def:curvature}
Consider a finite ground set $\mathcal{V}$ and a monotone submodular set function $f:2^\mathcal{V}\mapsto\mathbb{R}$ such that (without loss of generality) for each element $v \in \mathcal{V}$, it is $f(v)\neq 0$.  The curvature of $f$ is defined as $$\kappa_f=1-\min_{v\in\mathcal{V}}\frac{f(\mathcal{V})-f(\mathcal{V}\setminus\{v\})}{f(v)}.$$

\vspace{-3mm}
\hfill $\blacktriangle$
\end{mydef}

In words, the curvature of a monotone submodular function  $f:2^\mathcal{V}\mapsto\mathbb{R}$ measures how far $f$ is from modularity: in~particular, per Definition~\ref{def:sub} of submodularity, it follows that curvature takes values $0\leq\kappa_f\leq 1$, and
\begin{itemize}[leftmargin=*]
\item $\kappa_f=0$ if and only if for all elements $v\in\mathcal{V}$, it is $f(\mathcal{V})-f(\mathcal{V}\setminus\{v\})=f(v)$, that is, $f$ is modular.
\item $\kappa_f=1$ if and only if there exist an element $v\in\mathcal{V}$ such that $f(\mathcal{V})=f(\mathcal{V}\setminus\{v\})$, that is, in the presence of $\mathcal{V}\setminus\{v\}$, $v$ loses all its contribution to the overall value of $f(\mathcal{V})$.
\end{itemize} 

An example of a monotone submodular function with zero total curvature is the trace of the controllability matrix, which captures the control effort to drive the system in the state space~\cite{bullo2014}. A function with curvature $1$ is the matroid rank function~\cite{iyer2013curvature}. At the same time, many practically interesting functions have curvature strictly smaller than $1$, such as the concave over modular functions~\cite[Section~2.1]{iyer2013curvature}, and the $\log\det$ of positive-definite matrices~\cite{sviridenko2015optimal,sharma2015greedy}, which are used in applications such as speech processing~\cite{lin2011class}, computer vision~\cite{jegelka2011submodularity},
feature selection for Gaussian process regression~\cite{krause2008near}, and sensor scheduling for target tracking~\cite{tzoumas2016scheduling}.

The notion of total curvature has served to tighten bounds for several submodular maximization problems, e.g.,  for the non-resilient optimization problem 
$\max_{\mathcal{A}\subseteq \mathcal{V}, |\mathcal{A}|\leq \alpha}f(\mathcal{A})$ the approximation bound of the greedy is tightened from the value $1-1/e$ to $\frac{1}{\kappa_f}\left(1-e^{-\kappa_f}\right)$~\cite{conforti1984curvature,vondrak2010submodularity,iyer2013curvature}.
Nonetheless, for resilient submodular maximization problems, such as Problem~\ref{pr:robust_sub_max}, the role of the curvature has not been addressed yet.  We provide the first results to this end in the next section.

\subsection{Algorithm for resilient submodular maximization}\label{subsec:main_theorem}

\begin{algorithm}[t]
\caption{Algorithm for Problem~\ref{pr:robust_sub_max}.}
\begin{algorithmic}[1]
\REQUIRE Per Problem~\ref{pr:robust_sub_max}, three are the inputs to Algorithm~\ref{alg:rob_sub_max}:
\begin{itemize}[leftmargin=*]
\item finite ground set $\mathcal{V}=\{v_1, v_2, \ldots, v_m\}$, where $m=|\mathcal{V}|$;
\item submodular and monotone set function $f:2^\mathcal{V} \mapsto \mathbb{R}$ such that $f$ is non-negative and $f(\emptyset)=0$;
\item and two integers $\alpha$ and $\beta$ such that $0\leq\beta \leq \alpha\leq |\mathcal{V}|$.
\end{itemize}
\ENSURE  Set $\mathcal{A}_\textsc{Res}$, with properties per  Theorem~\ref{th:alg_rob_sub_max_performance}.
\STATE $\mathcal{A}_1\leftarrow\emptyset, \mathcal{A}_2\leftarrow\emptyset$
\STATE Sort elements in $\mathcal{V}$ such that $\mathcal{V}=\{v_1', v_2', \ldots, v_m'\}$ and $f(v_1')\geq f(v_2')\geq \ldots \geq f(v_m')$
\STATE $\mathcal{A}_1\leftarrow\{v_1', v_2', \ldots, v_\beta'\}$
\WHILE {$|\mathcal{A}_2| < \alpha-\beta$}  
\STATE $x\in \arg\max_{y \in \mathcal{V}\setminus (\mathcal{A}_1\cup\mathcal{A}_2)}f(y|\mathcal{A}_2)$

\STATE $\mathcal{A}_2 \leftarrow \{x\}\cup \mathcal{A}_2$

\ENDWHILE
\STATE $\mathcal{A}_\textsc{Res}\leftarrow \mathcal{A}_1 \cup \mathcal{A}_2$
\end{algorithmic}\label{alg:rob_sub_max}
\end{algorithm}

We exploit the curvature Definition~\ref{def:curvature} to provide for the resilient maximization Problem~\ref{pr:robust_sub_max}  Algorithm~\ref{alg:rob_sub_max}.   Algorithm~\ref{alg:rob_sub_max} returns a solution for Problem~\ref{pr:robust_sub_max}, denoted by $\mathcal{A}_\textsc{Res}$, in two steps: first, in lines 1-3 Algorithm~\ref{alg:rob_sub_max} selects a set $\mathcal{A}_1$, which is composed of $\beta$ elements from the ground set~$\mathcal{V}$. Specifically, per line 2, each element $v\in\mathcal{A}_1$ is such that for all elements $v'\in \mathcal{V}\setminus\mathcal{A}_1$, it is $f(v)\geq f(v')$.  Second, \mbox{in lines 4-8}, Algorithm~\ref{alg:rob_sub_max} selects greedily from the set $\mathcal{V}\setminus \mathcal{A}_1$ a set $\mathcal{A}_2$, which is composed of $\alpha-\beta$ elements, and then, in line~8, Algorithm~\ref{alg:rob_sub_max} returns set $\mathcal{A}_\textsc{Res}$ as the union of \mbox{$\mathcal{A}_1$ and $\mathcal{A}_2$}.

Algorithm~\ref{th:alg_rob_sub_max_performance}'s performance is quantified in Theorem~\ref{th:alg_rob_sub_max_performance}, whose proof can be found in the appendix. {The intuition behind Algorithm~\ref{alg:rob_sub_max} is discussed in Section~\ref{subsec:intuition}.}

\begin{mytheorem}\label{th:alg_rob_sub_max_performance}
Per Problem~\ref{pr:robust_sub_max}, let
\begin{itemize}[leftmargin=*]
\item the real number $f^\star$ equal to the (optimal) value of Problem~\ref{pr:robust_sub_max}, i.e., $f^\star={\max}_{\mathcal{A} \subseteq \mathcal{V}, |\mathcal{A}| \leq \alpha}\; {\min}_{\mathcal{B} \subseteq \mathcal{A}, |\mathcal{B}| \leq \beta}\;f(\mathcal{A}\setminus \mathcal{B});$ 
\item and for any set $\mathcal{A}\subseteq \mathcal{V}$, the set $\mathcal{B}^\star(\mathcal{A})$ equal to the optimal set-removal of $\beta$ elements from $\mathcal{A}$, i.e., $\mathcal{B}^\star(\mathcal{A})\in\min_{\mathcal{B} \subseteq \mathcal{A}, |\mathcal{B}| \leq \beta}f(\mathcal{A}\setminus \mathcal{B}).$ 
\end{itemize}

The following two hold on Algorithm~\ref{alg:rob_sub_max}'s performance:

\begin{enumerate}[leftmargin=*]
\item Algorithm~\ref{alg:rob_sub_max} returns a set $\mathcal{A}_\textsc{Res} \subseteq \mathcal{V}$ such that $|\mathcal{A}_\textsc{Res}|\leq \alpha$, and
\begin{align*}
&f(\mathcal{A}_\textsc{Res}\setminus \mathcal{B}^\star(\mathcal{A}_\textsc{Res}))\geq \\ &\;\;\;\max\left(1-\kappa_f, \frac{1}{\beta+1}\right)\frac{1}{\kappa_f}\left(1-e^{-\kappa_f}\right)f^\star,
\end{align*}
and in particular, for $\kappa_f=0$, $f(\mathcal{A}_\textsc{Res}\setminus \mathcal{B}^\star(\mathcal{A}_\textsc{Res}))= f^\star$.

\item Algorithm~\ref{alg:rob_sub_max} runs in $O(|\mathcal{V}|(\alpha-\beta))$ evaluations of $f$. \hfill $\blacktriangle$
\end{enumerate}
\end{mytheorem}

\begin{figure}[t]
\begin{center}
\begin{tikzpicture}
\begin{axis}[
    axis lines = left,
    xlabel = $\kappa_f$,
    ylabel = {$g(\kappa_f)$},
    ymajorgrids=true,
    grid style=dashed,
    legend style={at={(1.01,1)}},
    ymin=0, ymax=1,
]
\addplot [
    domain=0.01:1, 
    samples=30, 
    color=blue,
    mark=square,
    ]
    {1.05*(1-x)/x*(1-exp(-x))+0.005};
\addlegendentry{$g(\kappa_f)=\frac{1-\kappa_f}{\kappa_f}\left(1-e^{-\kappa_f}\right)$}
\end{axis}
\end{tikzpicture}
\caption{Plot of $g(\kappa_f)$ versus curvature $\kappa_f$ of a monotone submodular function $f$. By definition, the curvature $\kappa_f$ of a monotone submodular function $f$ takes values  between $0$ and $1$. $g(\kappa_f)$ increases from $0$ to $0.5$ as $\kappa_f$ decreases from $1$ to $0$.
}\label{fig:bounds}
\end{center}
\end{figure}
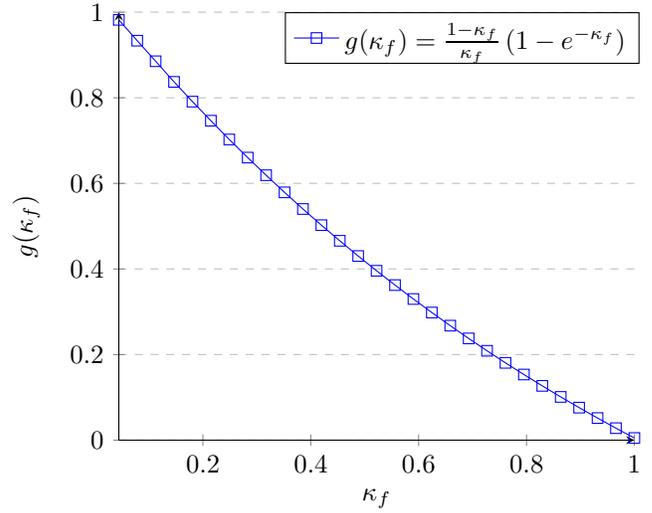


\begin{myremark}
Given any finite ground set $\mathcal{V}$ and monotone submodular function $f:2^\mathcal{V}\mapsto\mathbb{R}$ in the resilient maximization Problem~\ref{pr:robust_sub_max}, the approximation bound of Algorithm~\ref{alg:rob_sub_max} depends on the curvature value which is computed with $O(|\mathcal{V}|)$ evaluations of $f$ according to its Definition~\ref{def:curvature}.  \hfill $\blacktriangle$
\end{myremark}

\begin{myremark}
Given any finite ground set $\mathcal{V}$, finite number of failures or attacks $\beta$ and monotone submodular function $f:2^\mathcal{V}\mapsto\mathbb{R}$ in the resilient maximization Problem~\ref{pr:robust_sub_max}, the approximation bound of Algorithm~\ref{alg:rob_sub_max} is non-zero, since for finite $\beta$ it is $\max\left(1-\kappa_f, \frac{1}{\beta+1}\right)\geq \frac{1}{\beta+1}>0$, and for all $0\leq\kappa_f\leq 1$, it is $\frac{1}{\kappa_f}\left(1-e^{-\kappa_f}\right)\geq 1-1/e\simeq0.63>0$.  \hfill $\blacktriangle$
\end{myremark}

\begin{myremark}\label{rem:linear_approx}
Per Theorem~\ref{th:alg_rob_sub_max_performance},  Algorithm~\ref{alg:rob_sub_max}'s approximation performance bound is in the worst case equal to $\frac{1-\kappa_f}{\kappa_f}\left(1-e^{-\kappa_f}\right)$, which is plotted in Fig.~\ref{fig:bounds}, and is approximately equal to $-\kappa_f+1$. \hfill $\blacktriangle$
\end{myremark}

\begin{myremark}
For zero curvature, Algorithm~\ref{alg:rob_sub_max} solves Problem~\ref{pr:robust_sub_max} exactly. For non-zero curvature, in which case Problem~\ref{pr:robust_sub_max} is NP-hard~\cite{orlin2015robust}, \mbox{Algorithm~\ref{alg:rob_sub_max}'s} approximation bound is in the worst-case equal to $\frac{1-\kappa_f}{\kappa_f}\left(1-e^{-\kappa_f}\right)$, which tends to~$1$ as $\kappa_f\longrightarrow 0$. Overall, Algorithm~\ref{alg:rob_sub_max}'s curvature-dependent approximation bound makes a first step towards separating the class of monotone submodular functions into the functions for which the resilient maximization Problem~\ref{pr:robust_sub_max} can be approximated well (low curvature functions), and the functions for which it cannot (high curvature functions).  The~reason is that Algorithm~\ref{alg:rob_sub_max}'s approximation bound increases as the curvature decreases, and for zero curvature it becomes $1$ ---i.e., for zero curvature, Algorithm~\ref{alg:rob_sub_max} solves Problem~\ref{pr:robust_sub_max} exactly. 
This~role of the curvature in Problem~\ref{pr:robust_sub_max} is similar to the role that the curvature has played in the non-resilient variant of Problem~\ref{pr:robust_sub_max}, i.e., the optimization problem $\max_{\mathcal{A}\subseteq \mathcal{V}, |\mathcal{A}|\leq \alpha} \; f(\mathcal{A})$, where it has been used to separate the class of submodular functions into the functions for which $\max_{\mathcal{A}\subseteq \mathcal{V}, |\mathcal{A}|\leq \alpha} \; f(\mathcal{A})$  can be approximated well (low curvature functions), and the functions for which it cannot (high curvature functions)~\cite{conforti1984curvature,iyer2013curvature}.  Hence, Theorem~\ref{pr:robust_sub_max} also supports the intuition that the resilient maximization Problem~\ref{pr:robust_sub_max} is easier when the non-resilient variant $\max_{\mathcal{A}\subseteq \mathcal{V}, |\mathcal{A}|\leq \alpha} \; f(\mathcal{A})$ is easy, and it is harder when $\max_{\mathcal{A}\subseteq \mathcal{V}, |\mathcal{A}|\leq \alpha} \; f(\mathcal{A})$ is hard.

\hfill $\blacktriangle$ 
\end{myremark}

Theorem~\ref{th:alg_rob_sub_max_performance} implies for Algorithm~\ref{alg:rob_sub_max}'s performance:
\begin{itemize}[leftmargin=*]
\item Algorithm~\ref{alg:rob_sub_max} is the first scalable algorithm for Problem~\ref{pr:robust_sub_max} that is valid for any number of failures or attacks $\beta$, and any number of selections for resiliency $\alpha$.  In particular, the previously proposed algorithms for Problem~\ref{alg:rob_sub_max} in~\cite{krause2008robust} and~\cite{Orlin2016} are such that: the~algorithm in~\cite{krause2008robust} runs in exponential time in $\beta$, and quadratic in the cardinality of the ground set $\mathcal{V}$ and, as a result, has limited applicability in large-scale settings.  The~algorithm in~\cite{Orlin2016} is valid only for $\beta\leq \sqrt{2\alpha}$ and, as a result, has limited applicability in decision and control settings where the number of failures and attacks $\beta$ can be up to the number of placed sensors, actuators, etc. $\alpha$. For~example, the inequality $\beta\leq \sqrt{2\alpha}$ is violated when among $100$ placed sensors, $20$ may fail.

\item Algorithm~\ref{alg:rob_sub_max} is the first scalable algorithm for Problem~\ref{pr:robust_sub_max} that for non-zero curvature values $\kappa_f\leq 0.71$, and any number of failures and attacks $\beta$, guarantees approximation performance more than at least $29\%$ the optimal.
\mbox{In particular}, the previously proposed algorithms for Problem~\ref{alg:rob_sub_max} in~\cite{krause2008robust} and~\cite{Orlin2016} are such that: the algorithm in~\cite{krause2008robust} runs in exponential time in $\beta$, and quadratic in the cardinality of the ground set $\mathcal{V}$ and, as a result, has limited applicability in large-scale settings.  
The algorithm in~\cite{Orlin2016}, under the constraint $\beta\leq \sqrt{2\alpha}$, and for $\alpha, \beta \longrightarrow +\infty$, guarantees an approximate value up to at least $29\%$ the optimal. In~particular, for $\alpha, \beta < +\infty$, its approximation performance is less than at least $29\%$ the optimal.
\end{itemize}

An example of a central problem in machine learning, with applications, e.g., in sensor placement, where Algorithm~\ref{alg:rob_sub_max} guarantees almost exact approximation performance  ($\simeq 100\%$ the optimal) is that of Gaussian process regression for Gaussian processes with RBF kernels~\cite{krause2008near,bishop2006pattern}. The reason is that in this class of problems the objective function is the entropy of the selected measurements, which for Gaussian processes with RBF kernels was shown recently to have curvature values close to zero~\cite[Theorem~5]{sharma2015greedy}.

\subsection{Intuition behind curvature-dependent algorithm for resilient maximization}\label{subsec:intuition}

We explain the intuition behind Algorithm~\ref{alg:rob_sub_max}  for Problem~\ref{pr:robust_sub_max}, and give also an illustrative example. 
To this end, we focus only on the NP-hard case where the curvature of Problem~\ref{pr:robust_sub_max}'s objective function is non-zero~\cite[Lemma~3]{orlin2015robust}.  In addition, we use the notation introduced in the statements of Algorithm~\ref{alg:rob_sub_max} and of Theorem~\ref{th:alg_rob_sub_max_performance}.

We explain how Algorithm~\ref{alg:rob_sub_max} works first for the case where the optimal set-removal $\mathcal{B}^\star(\mathcal{A}_\textsc{Res})$ is equal to the set $\mathcal{A}_1$,  and then for the case where $\mathcal{B}^\star(\mathcal{A}_\textsc{Res})\neq \mathcal{A}_1$.  Notably, the case $\mathcal{B}^\star(\mathcal{A}_\textsc{Res})=\mathcal{A}_1$ is possible since in lines~1-3 Algorithm~\ref{alg:rob_sub_max} selects a set $\mathcal{A}_1$ such that $|\mathcal{A}_1|=\beta$, and per Problem~\ref{pr:robust_sub_max}, the~number {$\beta$ is the number of possible removals.}

\setcounter{paragraph}{0}
\paragraph{Intuition behind Algorithm~\ref{alg:rob_sub_max} for $\mathcal{B}^\star(\mathcal{A}_\textsc{Res})=\mathcal{A}_1$} 
The two parts of Algorithm~\ref{alg:rob_sub_max} operate in tandem to guarantee
\begin{equation}\label{bound:1}
f(\mathcal{A}_\textsc{Res}\setminus \mathcal{B}^\star(\mathcal{A}_\textsc{Res}))\geq \frac{1}{\kappa_f}\left(1-e^{-\kappa_f}\right) f^\star.
\end{equation}
This happens as follows: $\mathcal{B}^\star(\mathcal{A}_\textsc{Res})=\mathcal{A}_1$ implies that $\mathcal{A}_\textsc{Res}\setminus\mathcal{B}^\star(\mathcal{A}_\textsc{Res})=\mathcal{A}_2$, since $\mathcal{A}_\textsc{Res}=\mathcal{A}_1\cup\mathcal{A}_2$, per line~8 of Algorithm~\ref{alg:rob_sub_max}.  Therefore,  $f(\mathcal{A}_\textsc{Res}\setminus\mathcal{B}^\star(\mathcal{A}_\textsc{Res}))=f(\mathcal{A}_2)$.  But~in lines 4-7 of Algorithm~\ref{alg:rob_sub_max}, $\mathcal{A}_2$ is selected greedily and, as a result, using~\cite[Theorem~5.4]{conforti1984curvature} we have $$f(\mathcal{A}_2)\geq \frac{1}{\kappa_f}\left(1-e^{-\kappa_f}\right)\max_{\mathcal{A} \subseteq \mathcal{V}\setminus \mathcal{A}_1, |\mathcal{A}| \leq \alpha-\beta} f(\mathcal{A}).$$  Finally, it is $\max_{\mathcal{A} \subseteq \mathcal{V}\setminus \mathcal{A}_1, |\mathcal{A}| \leq \alpha-\beta}f(\mathcal{A})\geq f^\star,$ since the left-hand-side of this inequality is the maximum value one can achieve for $f$ by choosing $\alpha-\beta$ elements from the ground set $\mathcal{V}$ \emph{knowing} that $\mathcal{A}_1$ has been removed from $\mathcal{V}$, whereas the right-hand-side is the maximum value one can achieve for $f$ by choosing $\alpha$ elements from $\mathcal{V}$ \emph{not knowing} which $\beta$ of them will be optimally removed; a mathematical version of the latter proof can be found in~\cite[Lemma~2]{orlin2015robust}.
 
\paragraph{Intuition behind Algorithm~\ref{alg:rob_sub_max} for $\mathcal{B}^\star(\mathcal{A}_\textsc{Res})\neq\mathcal{A}_1$}
The two parts of Algorithm~\ref{alg:rob_sub_max} operate in tandem to guarantee,
\begin{align}\label{bound:2}
&f(\mathcal{A}_\textsc{Res}\setminus \mathcal{B}^\star(\mathcal{A}_\textsc{Res}))\geq \nonumber\\ &\;\;\;\max\left(1-\kappa_f, \frac{1}{\beta+1}\right)\frac{1}{\kappa_f}\left(1-e^{-\kappa_f}\right)f^\star,
\end{align}
This happens as follows: $\mathcal{B}^\star(\mathcal{A}_\textsc{Res})\neq\mathcal{A}_1$ implies, along with $|\mathcal{B}^\star(\mathcal{A}_\textsc{Res})|=|\mathcal{A}_1|=\beta$, that if $\mu$ elements in $\mathcal{A}_1$ are \emph{not} included in $\mathcal{B}^\star(\mathcal{A}_\textsc{Res})$, exactly $\mu$ elements in $\mathcal{A}_2$ are included in $\mathcal{B}^\star(\mathcal{A}_\textsc{Res})$.  Therefore, $f(\mathcal{A}_\textsc{Res}\setminus \mathcal{B}^\star(\mathcal{A}_\textsc{Res}))$ can take a bounded value similar to the one in ineq.~\eqref{bound:1} only if the $\mu$ elements in $\mathcal{A}_1$ that are \emph{not} included in $\mathcal{B}^\star(\mathcal{A}_\textsc{Res})$ can \emph{compensate} for the $\mu$ elements in $\mathcal{A}_2$ that are included in $\mathcal{B}^\star(\mathcal{A}_\textsc{Res})$.  \mbox{For this reason}, in line~2, Algorithm~\ref{alg:rob_sub_max} chooses the elements in $\mathcal{A}_1$ so that they a have higher value than those in $\mathcal{V}\setminus\mathcal{A}_1$: \mbox{In particular}, using the fact that the elements in $\mathcal{A}_1$ have a higher value than those in $\mathcal{V}\setminus\mathcal{A}_1$, and the definition of the curvature $\kappa_f$, in the proof of Theorem~\ref{th:alg_rob_sub_max_performance} we bound how much value the elements in $\mathcal{A}_1$ that are \emph{not} included in $\mathcal{B}^\star(\mathcal{A}_\textsc{Res})$ can compensate for the value of the elements in $\mathcal{A}_2$ that are included in $\mathcal{B}^\star(\mathcal{A}_\textsc{Res})$, and conclude ineq.~\eqref{bound:2}.

Overall, in the worst-case Algorithm~\ref{alg:rob_sub_max} guarantees for the resilient maximization Problem~\ref{pr:robust_sub_max} the approximation performance bound in ineq.~\eqref{bound:2}, as stated in Theorem~\ref{th:alg_rob_sub_max_performance}, since for all values of the curvature $\kappa_f$,  the bound in ineq.~\eqref{bound:2} is smaller than the bound in ineq.~\eqref{bound:1}, since $\max\left[1-\kappa_f, 1/(\beta+1)\right]\leq 1$, and we do not know a priori which of the two preceding bounds holds at a problem instance. 

\begin{myexample}\label{ex:illustration_of_algorithm}
We use an instance of Problem~\ref{pr:robust_sub_max} to illustrate how Algorithm~\ref{alg:rob_sub_max} finds an approximate solution to Problem~\ref{pr:robust_sub_max}, as well as, how it performs. We consider the following instance: let $\alpha=2$, $\beta=1$, $\mathcal{V}=\{v_1,v_2,v_3\}$, and $f$ such that $f(\emptyset)=0$, $f(v_3)>0$, $f(v_1)=f(v_3)+1$, $f(v_2)=f(v_3)+0.5$, $f(\{v_1\}\cup \{v_2\})=f(v_3)+1$, $f(\{v_1\}\cup \{v_3\})=f(\mathcal{V})=2f(v_3)+1$, and $f(\{v_2\}\cup \{v_3\})=2f(v_3)+0.5$.  

For the aforementioned instance, the curvature  is $\kappa_f=1$ and, as a result, Algorithm~\ref{alg:rob_sub_max} is guaranteed to return a set $\mathcal{A}_\textsc{Res}$ such that the approximation ratio is either at least $1-1/e$, per bound in ineq.~\eqref{bound:1}, or at least $(1-1/e)/2$, per bound in ineq.~\eqref{bound:2}. 

Algorithm~\ref{alg:rob_sub_max} selects $\mathcal{A}_\textsc{Res}=\{v_1,v_2\}$, which in this example is the exact solution to Problem~\ref{pr:robust_sub_max}.
In~particular, in lines 2-3, Algorithm~\ref{alg:rob_sub_max} selects $\mathcal{A}_1=\{v_1\}$, and in lines 4-7, it selects $\mathcal{A}_2=\{v_2\}$. \mbox{The optimal set-removal is $\mathcal{B}^\star(\mathcal{A}_\textsc{Res})=\{v_1\}$.}

The reason that Algorithm~\ref{alg:rob_sub_max} performs optimally in this example, which is not expected by Theorem~\ref{th:alg_rob_sub_max_performance} since it has the worst curvature value $\kappa_f=1$, is as follows: \mbox{In lines~1-3} Algorithm~\ref{alg:rob_sub_max} selects $\mathcal{A}_1=\{v_1\}$, which for $\mathcal{A}_\textsc{Res}=\{v_1,v_2\}$ is the element that  will be included in the optimal removal $\mathcal{B}^\star(\mathcal{A}_\textsc{Res})$. That is, in this example $\mathcal{B}^\star(\mathcal{A}_\textsc{Res})=\mathcal{A}_1$, which implies that the approximation performance of Algorithm~\ref{alg:rob_sub_max} is bounded by $1-1/e$, as in ineq.~\eqref{bound:1}. This is the first important observation towards explaining the optimal performance of Algorithm~\ref{alg:rob_sub_max} in this example; the second and final necessary observation is as follows: In lines~4-7, Algorithm~\ref{alg:rob_sub_max} selects the best element in $\mathcal{V}$ assuming that the element in $\mathcal{A}_1$ will be included in $\mathcal{B}^\star(\mathcal{A}_\textsc{Res})$, i.e., removed from $\mathcal{A}_\textsc{Res}$, since it selects greedily from $\mathcal{V}\setminus\mathcal{A}_1$. In contrast, if in lines~4-7 Algorithm~\ref{alg:rob_sub_max} would select greedily from $\mathcal{V}$ without taking into account that the element in $\mathcal{A}_1$ is going to be included in $\mathcal{B}^\star(\mathcal{A}_\textsc{Res})$, then it would select $\mathcal{A}_\textsc{Res}=\{v_1,v_3\}$ which is suboptimal for Problem~\ref{alg:rob_sub_max}.
\hfill $\blacktriangle$
\end{myexample}

\section{Simulations} \label{sec:experiments}

We empirically test Algorithm~\ref{alg:rob_sub_max}'s approximation performance for Problem~\ref{pr:robust_sub_max} against an exact, brute force algorithm. As a test function, with consider one of the following form: Given a finite ground set $\mathcal{V}$, and $|\mathcal{V}|$ positive semi-definite matrices $D_1, D_2, \ldots, D_{|\mathcal{V}|}$, for any set $\mathcal{A}\subseteq \mathcal{V}$, let the set function $f(\mathcal{A})=\log\det(\sum_{i\in\mathcal{A}}D_i+I)$, where $I$ is the identity matrix. Functions of this form appear in applications such as sensor selection for Gaussian process regression~\cite{krause2008near}, and actuator placement for bounded control effort~\cite{summers2014submodularity,tzoumas2014minimal}. \mbox{To run our simulations}, and be able to compute the exact value to Problem~\ref{pr:robust_sub_max}, we select small sizes of $|\mathcal{V}|$ from $8$ to $15$. In addition, we fix the number of selections for resiliency $\alpha$ to $7$, vary the number of attacks/failures $\beta$ from $1$ to $6$, and for each of the aforementioned cases, generate $10$ random instances of the matrices $D_1, D_2, \ldots, D_{|\mathcal{V}|}$ of size $20\times 20$. 

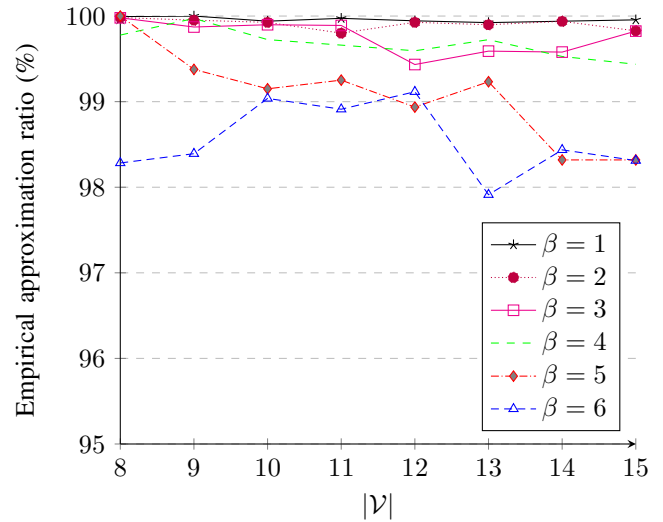
\begin{figure}[t]
\begin{center}
\begin{tikzpicture}
\begin{axis}[
    axis lines = left,
    ymin=95, ymax=100,
    xlabel = $|\mathcal{V}|$,
    ylabel = {\text{Empirical approximation ratio (\%)}},
    legend pos=south east,
    ymajorgrids=true,
    grid style=dashed,
]

\addplot[
    color=black,
    mark=cross,
    style={solid},mark=star,
    ]
    coordinates { 
    (8,99.9938969591931)(9,100)(10,99.9403036983340
    )(11,99.9734405209467)(12,99.9445055851978)(13,99.9233569269559)(14,99.9399827533769)(15,99.9556339382374)
    };
\addlegendentry{$\beta=1$}
\addplot[
    color=purple,
    mark=cross, 
    style={densely  dotted}, mark=otimes*
    ]
    coordinates {  
    (8,99.9823778779938)(9,99.9518916555192)(10,99.9265206910867
    )(11,99.8002556085342)(12,99.9269852829128)(13,99.8986100587955)(14,99.9385511168707)(15,99.8272566728878)
    };
\addlegendentry{$\beta=2$}
\addplot[
    color=magenta,
    mark=cross,
    style={solid},mark=square
    ]
    coordinates {  
    (8,99.9795756218319)(9,99.8720812573521)(10,99.8986950608159
    )(11,99.8900252193008)(12,99.4337328361950)(13,99.5902722668714)(14,99.5791309111200)(15,99.8243471690982)
    };
\addlegendentry{$\beta=3$}
\addplot[
    color=green,
    style={dashed},
    ]
    coordinates {
    (8,99.7776213681899)(9,99.9721189548437)(10,99.7248207194093
    )(11,99.6604324797294)(12,99.5942004293258)(13,99.7243480330286)(14,99.5276336910729)(15,99.4367575405263)
    };
\addlegendentry{$\beta=4$}
\addplot[
    color=red,
    mark=cross,densely dashdotted,every mark/.append style={solid, fill=gray},mark=diamond*
    ]
    coordinates {
    (8,100)(9,99.3767595556905)(10,99.1500311161459
    )(11,99.2514709445560)(12,98.9363423110778)(13,99.2326335209572)(14,98.3191500252632)(15,98.3191500252632)
    };
\addlegendentry{$\beta=5$}
\addplot[
    color=blue,densely dashed, every mark/.append style={solid},mark=triangle
    ]
    coordinates {
    (8,98.2836281708944)(9,98.3917870409489)(10,99.0376637971792)(11,98.9115422065165)(12,99.1157002951941)(13,97.9130884056088)(14,98.4367108676951)(15,98.3117129289644)
    };
\addlegendentry{$\beta=6$}
\end{axis}
\end{tikzpicture}
\caption{Empirical approximation performance of Algorithm~\ref{alg:rob_sub_max}. For details, see Section~\ref{sec:experiments}.
}\label{fig:sim}
\end{center}
\end{figure}

Our simulations are summarized in Fig.~\ref{fig:sim}, where Algorithm~\ref{alg:rob_sub_max} is seen to perform close to $98\%$ the optimal, and its approximation performance to degrade as $\beta$ increases up to~$\alpha$, which is equal to $7$. Notably, for all generated instances of~$f$, $f$'s curvature takes values larger than $0.9$, for which, Algorithm~\ref{alg:rob_sub_max}'s  worst-case theoretical approximation bound in Theorem~\ref{th:alg_rob_sub_max_performance} is at least $14\%$, since for the maximum value of~$\beta$, which is equal to $6$ in this simulation example, Algorithm~\ref{alg:rob_sub_max}'s approximation bound per Theorem~\ref{th:alg_rob_sub_max_performance} is  $ \max[1/(\beta+1), 1-\kappa_f]=1/(\beta+1)=0.14$. Overall, Algorithm~\ref{alg:rob_sub_max}'s approximation performance in Fig.~\ref{fig:sim} is in accordance with the empirical observation that greedy-type algorithms for submodular maximization often outperform in practice their worst-case theoretical approximation guarantees~\cite{krause2008near,sharma2015greedy}.

\section{Concluding remarks \& Future work} \label{sec:con}

We focused on the resilient submodular maximization Problem~\ref{pr:robust_sub_max}, 
which is central in machine learning, optimization, and control, in applications such as feature selection for classifier training, and sensor scheduling for target tracking.
In particular, exploiting the notion of curvature, we provided the first scalable algorithm for Problem~\ref{pr:robust_sub_max}, Algorithm~\ref{alg:rob_sub_max}, which is valid for any number of attacks or failures, and which, for functions with low curvature, guarantees superior approximation performance.  In addition, for functions with zero curvature, Algorithm~\ref{alg:rob_sub_max} solves Problem~\ref{pr:robust_sub_max} exactly. 
Overall, Algorithm~\ref{alg:rob_sub_max}'s approximation bound makes a first step to characterize the curvature's effect on approximations for resilient submodular maximization problems, complementing that way the current knowledge on the curvature's effect on \emph{non-resilient} submodular maximization.
Future work will focus on Algorithm~\ref{alg:rob_sub_max}'s \mbox{extension to matroid constraints.}

\section*{Appendix: Proof of Theorem~\ref{th:alg_rob_sub_max_performance}}

\paragraph*{Notation} Given a set function $f:2^\mathcal{V}\mapsto \mathbb{R}$, for any sets $\mathcal{X}\subseteq \mathcal{V}$ and $\mathcal{X}'\subseteq \mathcal{V}$, the $f(\mathcal{X}|\mathcal{X}')$ denotes $f(\mathcal{X}\cup\mathcal{X}')-f(\mathcal{X}')$. The set $\mathcal{A}^\star$ denotes a solution to Problem~\ref{pr:robust_sub_max}, i.e., $\mathcal{A}^\star\in \arg{\max}_{\mathcal{A} \subseteq \mathcal{V}, |\mathcal{A}| \leq \alpha}\; {\min}_{\mathcal{B} \subseteq \mathcal{A}, |\mathcal{B}| \leq \beta}\;f(\mathcal{A}\setminus \mathcal{B})$.

\subsection{Proof of Algorithm~\ref{alg:rob_sub_max}'s approximation performance:} 
We complete the proof first for curvature value $\kappa_f$ equal to $0$, and then for~$\kappa_f\neq 0$.

\subsubsection{Proof of Algorithm~\ref{alg:rob_sub_max}'s approximation performance for $\kappa_f=0$} In this case, Algorithm~\ref{alg:rob_sub_max} solves Problem~\ref{pr:robust_sub_max} exactly. This follows from the two observations below:
\begin{itemize}[leftmargin=*]
\item For $\kappa_f=0$, Algorithm~\ref{alg:rob_sub_max} returns a set $\mathcal{A}_\textsc{Res}$ such that for all elements $v\in\mathcal{A}_\textsc{Res}$ and $v'\in\mathcal{V}\setminus\mathcal{A}_\textsc{Res}$, it is $f(v)\geq f(v')$.  The reasons are two: first,the set~$\mathcal{A}_1$ is such that for all elements $v\in\mathcal{A}_1$ and $v'\in\mathcal{V}\setminus\mathcal{A}_1$, $f(v)\geq f(v')$; and second~the set $\mathcal{A}_2$ is such that for all elements $v\in\mathcal{A}_2$ and $v'\in\mathcal{V}\setminus(\mathcal{A}_1\cup \mathcal{A}_2)$, it is $f(v)\geq f(v')$, since in line~5 of Algorithm~\ref{alg:rob_sub_max} it is  $f(y|\mathcal{A}_1\cup\mathcal{A}_2)=f(y)$, which is implied from the fact that $\kappa_f=0$ implies that $f$ is modular.

\item For $\kappa_f=0$, the set $\mathcal{A}^\star$ is also such that for all elements $v\in\mathcal{A}^\star$ and $v'\in\mathcal{V}\setminus\mathcal{A}^\star$, it is $f(v)\geq f(v')$.  To explain why, we first make two observations: first,~$\kappa_f=0$ implies that $f$ is modular, which in turn implies that for any sets $\mathcal{A}\subseteq \mathcal{V}$ and $\mathcal{B}\subseteq \mathcal{A}$, it is  $f(\mathcal{A}\setminus\mathcal{B})=\sum_{v\in \mathcal{A}\setminus\mathcal{B}}f(v)$; and second, Problem~\ref{pr:robust_sub_max} considers (without loss of generality) that $f$ is non-negative, which implies that for all $v \in \mathcal{V}$, it is $f(v)\geq 0$.  Therefore, for any $\mathcal{A}$, the optimal set-removal $\mathcal{B}^\star(\mathcal{A})$ contains $\beta$ elements in $\mathcal{A}$ such that for all elements $x\in \mathcal{B}^\star(\mathcal{A})$ and $y\in \mathcal{A}$, it is $f(x)\geq f(y)$.  Thus, $\mathcal{A}^\star$ maximizes $f(\mathcal{A}\setminus\mathcal{B}^\star(\mathcal{A}))$ if and only if for all elements $v\in\mathcal{A}$ and $v'\in\mathcal{V}\setminus\mathcal{A}$, it is $f(v)\geq f(v')$.
\end{itemize}

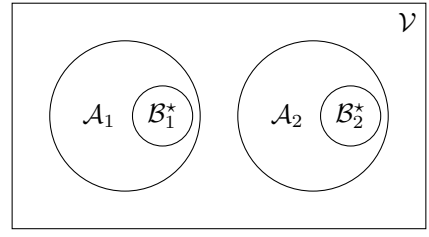
\begin{figure}[t]
\def \setAone{ (0,0) circle (1cm) }
\def \setBone{ (.5,0) circle (0.4cm)}
\def \setAtwo{ (2.5,0) circle (1cm) }
\def \setBtwo{ (3.0,0) circle (0.4cm)}
\def \myrectangle{ (-1.5, -1.5) rectangle (4, 1.5) }
\begin{center}
\begin{tikzpicture}
\draw \myrectangle node[below left]{$\mathcal{V}$};
\draw \setAone node[left]{$\mathcal{A}_1$};
\draw \setBone node[]{$\mathcal{B}_1^\star$};
\draw \setAtwo node[left]{$\mathcal{A}_2$};
\draw \setBtwo node[]{$\mathcal{B}_2^\star$};
\end{tikzpicture}
\end{center}
\caption{Venn diagram, where $\mathcal{A}_1, \mathcal{A}_2,\mathcal{B}_1^\star, \mathcal{B}_2^\star$ are as follows: Per Algorithm~\ref{alg:rob_sub_max}, $\mathcal{A}_1$  and $\mathcal{A}_2$ are such that $\mathcal{A}_\textsc{Res}=\mathcal{A}_1\cup \mathcal{A}_2$, and $\mathcal{A}_1\cap \mathcal{A}_2=\emptyset$. Also, 
$\mathcal{B}_1^\star$ and $\mathcal{B}_2^\star$ are such that  $\mathcal{B}_1^\star=\mathcal{B}^\star(\mathcal{A}_\textsc{Res})\cap\mathcal{A}_1$, and $\mathcal{B}_2^\star=\mathcal{B}^\star(\mathcal{A}_\textsc{Res})\cap\mathcal{A}_2$; by definition, $\mathcal{B}_1^\star\cap \mathcal{B}_2^\star=\emptyset$ and $\mathcal{B}^\star(\mathcal{A}_\textsc{Res})=\mathcal{B}_1^\star\cup \mathcal{B}_2^\star$.  
}\label{fig:venn_diagram_for_proof}
\end{figure}

\subsubsection{Proof of Algorithm~\ref{alg:rob_sub_max}'s approximation performance for $\kappa_f\neq 0$} 
We use the symbols $\mathcal{B}_1^\star$ and~$\mathcal{B}_2^\star$, defined in Fig.~\ref{fig:venn_diagram_for_proof}.
 
We complete the proof of Algorithm~\ref{alg:rob_sub_max}'s approximation bound for $\kappa_f\neq 0$ in two steps: First, we consider that  $\mathcal{B}^\star(\mathcal{A}_\textsc{Res})=\mathcal{A}_1$, i.e., $\mathcal{B}_2^\star=\emptyset$. Second, we consider that $\mathcal{B}^\star(\mathcal{A}_\textsc{Res})\neq\mathcal{A}_1$, i.e., $\mathcal{B}_2^\star\neq \emptyset$.

For the case $\mathcal{B}_2^\star= \emptyset$, we prove 
\begin{equation}\label{ineq:optimality_result_first_case}
f(\mathcal{A}_\textsc{Res}\setminus \mathcal{B}^\star(\mathcal{A}_\textsc{Res}))\geq \frac{1}{\kappa_f}\left(1-e^{-\kappa_f}\right) f(\mathcal{A}^\star\setminus \mathcal{B}^\star(\mathcal{A}^\star)).
\end{equation}  
In particular, we prove ineq.~\eqref{ineq:optimality_result_first_case} by making the following three consecutive observations:  first, $f(\mathcal{A}_\textsc{Res}\setminus \mathcal{B}^\star(\mathcal{A}_\textsc{Res}))=f(\mathcal{A}_2)$, since $\mathcal{B}^\star(\mathcal{A}_\textsc{Res})=\mathcal{A}_1$, because $\mathcal{B}_2^\star= \emptyset$.
Second,
\begin{equation}\label{ineq:lower_bound_of_f(A_2)}
 f(\mathcal{A}_2)\geq \frac{1}{\kappa_f}\left(1-e^{-\kappa_f}\right)\max_{\mathcal{A} \subseteq \mathcal{V}\setminus \mathcal{A}_1, |\mathcal{A}| \leq \alpha-\beta} f(\mathcal{A}),
\end{equation}
since in line~5 Algorithm~\ref{alg:rob_sub_max} constructs greedily the set $\mathcal{A}_2$  using elements from $\mathcal{V}\setminus \mathcal{A}_1$~\cite[Theorem 5.4]{conforti1984curvature}; and third,
\begin{equation}\label{ineq:bound_of_A_2_by_opt}
\max_{\mathcal{A} \subseteq \mathcal{V}\setminus \mathcal{A}_1, |\mathcal{A}| \leq \alpha-\beta}f(\mathcal{A})\geq f(\mathcal{A}^\star\setminus \mathcal{B}^\star(\mathcal{A}^\star)),
\end{equation} 
because of~\cite[Lemma~2]{orlin2015robust}, where we recall that  the set $\mathcal{A}^\star$ denotes an (optimal) solution to Problem~\ref{pr:robust_sub_max}.
The above three consecutive observations complete the proof of ineq.~\eqref{ineq:optimality_result_first_case}, and end our focus on the case where $\mathcal{B}_2^\star= \emptyset$.

For the case $\mathcal{B}_2^\star\neq \emptyset$, we prove first that
\begin{align}\label{ineq:optimality_result_second_case}
&f(\mathcal{A}_\textsc{Res}\setminus \mathcal{B}^\star(\mathcal{A}_\textsc{Res}))\geq \frac{1-\kappa_f}{\kappa_f}\left(1-e^{-\kappa_f}\right) \nonumber\\
&\hspace{54mm}f(\mathcal{A}^\star\setminus \mathcal{B}^\star(\mathcal{A}^\star)),
\end{align} 
and then we complete the proof of Algorithm~\ref{alg:rob_sub_max}'s approximation performance bound by also proving that
\begin{align}\label{ineq:optimality_result_second_case_2}
&f(\mathcal{A}_\textsc{Res}\setminus \mathcal{B}^\star(\mathcal{A}_\textsc{Res}))\geq\nonumber\\
&\;\;\; \frac{1}{\beta+1}\frac{1}{\kappa_f}\left(1-e^{-\kappa_f}\right) f(\mathcal{A}^\star\setminus \mathcal{B}^\star(\mathcal{A}^\star)).
\end{align}  
   
To prove ineq.~\eqref{ineq:optimality_result_second_case}, we use the following lemma.

\begin{mylemma}\label{lem:curvature}
Consider a finite ground set $\mathcal{V}$ and a monotone set function $f:2^\mathcal{V}\mapsto \mathbb{R}$ such that $f$ is non-negative and $f(\emptyset)=0$. For any set $\mathcal{A}\subseteq \mathcal{V}$,
\belowdisplayskip=-12pt\begin{equation*}
f(\mathcal{A})\geq (1-\kappa_f)\sum_{a \in \mathcal{A}}f(a).
\end{equation*}
\hfill $\blacktriangle$
\end{mylemma}
\vspace{2mm}
\textit{Proof of Lemma~\ref{lem:curvature}:} Let $\mathcal{A}=\{a_1,a_2,\ldots, a_{|{\cal A}|}\}$. We prove Lemma~\ref{lem:curvature} by proving the following two inequalities: 
\begin{align}
f(\mathcal{A})&\geq \sum_{i=1}^{|{\cal A}|} f(a_i|\mathcal{V}\setminus\{a_i\}),\label{ineq:aux_5}\\
\sum_{i=1}^{|{\cal A}|} f(a_i|\mathcal{V}\setminus\{a_i\})&\geq (1-\kappa_f)\sum_{i=1}^{|{\cal A}|} f(a_i)\label{ineq:aux_6}. 
\end{align} 

We begin with the proof of ineq.~\eqref{ineq:aux_5}: 
\begin{align}
f(\mathcal{A})&=f(\mathcal{A}|\emptyset)\label{ineq:aux_9}\\
&\geq f(\mathcal{A}|\mathcal{V}\setminus \mathcal{A})\label{ineq:aux_10}\\
&= \sum_{i=1}^{|{\cal A}|}f(a_i|\mathcal{V}\setminus\{a_i,a_{i+1},\ldots,a_{|{\cal A}|}\})\label{ineq:aux_11}\\
&\geq \sum_{i=1}^{|{\cal A}|}f(a_i|\mathcal{V}\setminus\{a_i\}),\label{ineq:aux_12}
\end{align}
where ineqs.~\eqref{ineq:aux_10} to~\eqref{ineq:aux_12} hold for the following reasons: ineq.~\eqref{ineq:aux_10} is implied by eq.~\eqref{ineq:aux_9} because $f$ is submodular and $\emptyset\subseteq \mathcal{V}\setminus \mathcal{A}$; eq.~\eqref{ineq:aux_11} holds since for any sets $\mathcal{X}\subseteq \mathcal{V}$ and $\mathcal{Y}\subseteq \mathcal{V}$ it is $f(\mathcal{X}|\mathcal{Y})=f(\mathcal{X}\cup \mathcal{Y})-f(\mathcal{Y})$, and it also  $\{a_1,a_2,\ldots, a_{|{\cal A}|}\}$ denotes the set $\mathcal{A}$; and ineq.~\eqref{ineq:aux_12} holds since $f$ is submodular and $\mathcal{V}\setminus\{a_i,a_{i+1},\ldots,a_{\mu}\} \subseteq \mathcal{V}\setminus\{a_i\}$.  These observations complete the proof of ineq.~\eqref{ineq:aux_5}.

We now prove ineq.~\eqref{ineq:aux_6} using the Definition~\ref{def:curvature} of $\kappa_f$, as follows: since $\kappa_f=1-\min_{v\in \mathcal{V}}\frac{f(v|\mathcal{V}\setminus\{v\})}{f(v)}$, it is implied that for all elements $v\in \mathcal{V}$ it is $ f(v|\mathcal{V}\setminus\{v\})\geq (1-\kappa_f)f(v)$.  Therefore, adding the latter inequality across all elements $a \in \calA$ completes the proof of ineq.~\eqref{ineq:aux_6}.
\hfill $\square$

In addition to the above lemma, to prove ineq.~\eqref{ineq:optimality_result_second_case} we~use the following two notations: first, let the set $\mathcal{A}_1^+$ denote the elements in $\mathcal{A}_1$ not included in the optimal set-removal $\mathcal{B}^\star(\mathcal{A}_\textsc{Res})$; notably, $\mathcal{A}_1^+$~is non-empty, since $\mathcal{B}^\star(\mathcal{A}_\textsc{Res})$ intersects with the set $\mathcal{A}_2$ (because $\mathcal{B}_2^\star\neq \emptyset$) and $|\mathcal{B}^\star(\mathcal{A}_\textsc{Res})|=|\mathcal{A}_1|$ and, as a result, at least one element in $\mathcal{A}_1$ is not included in $\mathcal{B}^\star(\mathcal{A}_\textsc{Res})$.  In addition, let the set $\mathcal{A}_2^+$ denote the elements in $\mathcal{A}_2$ not included in $\mathcal{B}^\star(\mathcal{A}_\textsc{Res})$.

We prove ineq.~\eqref{ineq:optimality_result_second_case} first using Lemma~\ref{lem:curvature} and then taking into account ineqs.~\eqref{ineq:lower_bound_of_f(A_2)} and~\eqref{ineq:bound_of_A_2_by_opt}, as follows:
\begin{align}
& f(\mathcal{A}_\textsc{Res}\setminus\mathcal{B}^\star(\mathcal{A}_\textsc{Res})) \nonumber\\
&\;\;\;= f(\mathcal{A}_1^+\cup \mathcal{A}_2^+)\label{ineq:aux_14}\\
&\;\;\;\geq (1-\kappa_f)\sum_{v \in \mathcal{A}_1^+\cup \mathcal{A}_2^+}f(v)\label{ineq:aux_15}\\
&\;\;\;\geq (1-\kappa_f)\left(\sum_{v \in \mathcal{A}_2\setminus \mathcal{A}_2^+}f(v)+\sum_{v \in \mathcal{A}_2^+}f(v)\right)\label{ineq:aux_16}\\
&\;\;\;\geq (1-\kappa_f)f[(\mathcal{A}_2\setminus \mathcal{A}_2^+)\cup \mathcal{A}_2^+]\label{ineq:aux_17}\\
&\;\;\;= (1-\kappa_f)f(\mathcal{A}_2)\label{ineq:aux_18},
\end{align}
where eq.~\eqref{ineq:aux_14} to~\eqref{ineq:aux_18} hold for the following reasons: eq.~\eqref{ineq:aux_14} follows from the definitions of the sets~$\mathcal{A}_1^+$ and $\mathcal{A}_2^+$; ineq.~\eqref{ineq:aux_15} follows from ineq.~\eqref{ineq:aux_14} due to Lemma~\ref{lem:curvature}; ineq.~\eqref{ineq:aux_16} follows from ineq.~\eqref{ineq:aux_15} because for all elements $v \in \mathcal{A}_1^+$ and all elements  $v' \in \mathcal{A}_2\setminus\mathcal{A}_2^+$ it is $f(v)\geq f(v')$ (note that due to the definitions of the sets~$\mathcal{A}_1^+$ and $\mathcal{A}_2^+$ it is $|\mathcal{A}_1^+|=|\mathcal{A}_2\setminus\mathcal{A}_2^+|$, that is, the number of non-removed elements in $\calA_1$ is equal to the number of removed elements in $\calA_2$);  finally, ineq.~\eqref{ineq:aux_17} follows from ineq.~\eqref{ineq:aux_16} because the set function $f$ is submodular and, as~a result, the~submodularity Definition~\ref{def:sub} implies that for any sets $\mathcal{A}\subseteq \mathcal{V}$ and $\mathcal{A}'\subseteq \mathcal{V}$, it is  $f(\mathcal{A})+f(\mathcal{A}')\geq f(\mathcal{A}\cup \mathcal{A}')$.
Now, ineq.~\eqref{ineq:optimality_result_second_case}  follows from ineq.~\eqref{ineq:aux_18} by taking into account ineqs.~\eqref{ineq:lower_bound_of_f(A_2)} and~\eqref{ineq:bound_of_A_2_by_opt}.

To complete the proof of Algorithm~\ref{alg:rob_sub_max}'s approximation performance bound it remains to prove ineq.~\eqref{ineq:optimality_result_second_case_2}.  To this end, we denote
\begin{align}
	\eta = \frac{f(\mathcal{B}_2^\star|\mathcal{A}_\textsc{Res}\setminus \mathcal{B}^\star(\mathcal{A}_\textsc{Res}))}{f(\mathcal{A}_2)}
\end{align}
Later in this proof, we prove that $0\leq \eta\leq 1$.  We prove ineq.~\eqref{ineq:optimality_result_second_case_2} by first observing that
\begin{equation}\label{ineq:aux_1}
f(\mathcal{A}_\textsc{Res}\setminus\mathcal{B}^\star(\mathcal{A}_\textsc{Res}))\geq\max\{f(\mathcal{A}_\textsc{Res}\setminus\mathcal{B}^\star(\mathcal{A}_\textsc{Res})),f(\mathcal{A}_1^+)\},
\end{equation}
and then proving the following three inequalities:
\begin{align}
f(\mathcal{A}_\textsc{Res}\setminus\mathcal{B}^\star(\mathcal{A}_\textsc{Res}))&\geq(1-\eta)f(\mathcal{A}_2)\label{ineq:aux_2},\\
f(\mathcal{A}_1^+)&\geq \eta \frac{1}{\beta}f(\mathcal{A}_2),\label{ineq:aux_3}\\
\max\{(1-\eta),\eta\frac{1}{\beta}\}&\geq \frac{1}{\beta+1}.\label{ineq:aux_4}
\end{align}

Specifically, if we substitute ineqs.~\eqref{ineq:aux_2},~\eqref{ineq:aux_3} and~\eqref{ineq:aux_4} to~\eqref{ineq:aux_1}, and take into account that $f(\mathcal{A}_2)\geq 0$, then
\begin{equation*}
f(\mathcal{A}_\textsc{Res}\setminus\mathcal{B}^\star(\mathcal{A}_\textsc{Res}))\geq \frac{1}{\beta+1}f(\mathcal{A}_2),
\end{equation*}
which implies ineq.~\eqref{ineq:optimality_result_second_case_2} after taking into account ineqs.~\eqref{ineq:lower_bound_of_f(A_2)}~and~\eqref{ineq:bound_of_A_2_by_opt}.

In the remaining paragraphs, we complete the proof of ineq.~\eqref{ineq:optimality_result_second_case_2} by proving $0\leq \eta\leq 1$, and ineqs.~\eqref{ineq:aux_2},~\eqref{ineq:aux_3} and~\eqref{ineq:aux_4}, respectively.

\paragraph{Proof of ineq.~$0\leq \eta\leq 1$} We first prove that $\eta\geq 0$, and then that $\eta\leq 1$: i)~$\eta\geq 0$, since by definition $\eta=f(\mathcal{B}_2^\star|\mathcal{A}_\textsc{Res}\setminus \mathcal{B}^\star(\mathcal{A}_\textsc{Res}))/f(\mathcal{A}_2)$, and $f$ is non-negative; and ii)~$\eta\leq 1$, since $f(\mathcal{A}_2)\geq f(\mathcal{B}^\star_2)$, due to monotonicity of $f$ and that $\mathcal{B}^\star_2 \subseteq \mathcal{A}_2$, and $f(\mathcal{B}^\star_2)\geq f(\mathcal{B}_2^\star|\mathcal{A}_\textsc{Res}\setminus \mathcal{B}^\star(\mathcal{A}_\textsc{Res}))$, due to submodularity of $f$ and that $\emptyset \subseteq \mathcal{A}_\textsc{Res}\setminus \mathcal{B}^\star(\mathcal{A}_\textsc{Res})$. 


\paragraph{Proof of ineq.~\eqref{ineq:aux_2}}  We complete the proof of ineq.~\eqref{ineq:aux_2} in two steps.  First, it can be verified that
\begin{align}\label{eq:aux_1}
& f(\mathcal{A}_\textsc{Res}\setminus\mathcal{B}^\star(\mathcal{A}_\textsc{Res}))=f(\mathcal{A}_2)-\nonumber\\ & f(\mathcal{B}^\star_2|\mathcal{A}_\textsc{Res}\setminus\mathcal{B}^\star(\mathcal{A}_\textsc{Res}))+f(\mathcal{A}_1|\mathcal{A}_2)-f(\mathcal{B}^\star_1|\mathcal{A}_\textsc{Res}\setminus\mathcal{B}^\star_1),
\end{align}
since for any $\mathcal{X}\subseteq \mathcal{V}$ and $\mathcal{Y}\subseteq \mathcal{V}$, $f(\mathcal{X}|\mathcal{Y})=f(\mathcal{X}\cup \mathcal{Y})-f(\mathcal{Y})$. Second,~\eqref{eq:aux_1} implies~\eqref{ineq:aux_2}, since $f(\mathcal{B}^\star_2|\mathcal{A}_\textsc{Res}\setminus\mathcal{B}^\star(\mathcal{A}_\textsc{Res}))=\eta f(\mathcal{A}_2)$, and $f(\mathcal{A}_1|\mathcal{A}_2)-f(\mathcal{B}^\star_1|\mathcal{A}_\textsc{Res}\setminus\mathcal{B}^\star_1)\geq 0$.
The latter is true due to the following two observations: i)~$f(\mathcal{A}_1|\mathcal{A}_2)\geq f(\mathcal{B}_1^\star|\mathcal{A}_2)$, since $f$ is monotone and $\mathcal{B}_1^\star \subseteq \mathcal{A}_1$; and ii)~$f(\mathcal{B}_1^\star|\mathcal{A}_2)\geq f(\mathcal{B}^\star_1|\mathcal{A}_\textsc{Res}\setminus\mathcal{B}^\star_1)$, since $f$ is submodular and $\mathcal{A}_2\subseteq \mathcal{A}_\textsc{Res}\setminus\mathcal{B}^\star_1$ (see also Fig.~\ref{fig:venn_diagram_for_proof}).

\paragraph{Proof of ineq.~\eqref{ineq:aux_3}} 

We use the following lemma.
 

\begin{mylemma}\label{lem:D3}
Consider any finite ground set $\mathcal{V}$, a monotone submodular function $f:2^\mathcal{V}\mapsto \mathbb{R}$ and non-empty sets $\calY, \calP\subseteq \calV$ such that for all elements $y \in \calY$ and all elements $p \in \calP$ it is $f(y)\geq f(p)$.  Then,
\belowdisplayskip=-12pt\begin{equation*}
f(\calP|\calY)\leq |\calP|f(\calY).
\end{equation*}
\hfill $\blacktriangle$
\end{mylemma}
\vspace{2mm}
\textit{Proof of Lemma~\ref{lem:D3}:} Consider any element $y \in \calY$ (such an element exists since Lemma~\ref{lem:D3} considers that $\calY$ is non-empty); then, 
\begin{align}
f(\calP|\calY)&= f(\calP\cup\calY)-f(\calY)\label{aux1:1}\\
&\leq f(\calP)+f(\calY)-f(\calY)\label{aux1:2}\\
&= f(\calP)\nonumber\\
&\leq \sum_{p\in\calP}f(p)\label{aux1:4}\\
&\leq |\calP| \max_{p\in\calP} f(p)\nonumber\\
&\leq |\calP|  f(y)\label{aux1:7}\\
&\leq |\calP| f(\calY),\label{aux1:5}
\end{align}
where eq.~\eqref{aux1:1} to ineq.~\eqref{aux1:5} hold for the following reasons: eq.~\eqref{aux1:1} holds since for any sets $\mathcal{X}\subseteq \mathcal{V}$ and $\mathcal{Y}\subseteq \mathcal{V}$, it is $f(\mathcal{X}|\mathcal{Y})=f(\mathcal{X}\cup \mathcal{Y})-f(\mathcal{Y})$; ineq.~\eqref{aux1:2} holds since $f$ is submodular and, as a result, the submodularity Definition~\ref{def:sub} implies that for any set $\calA\subseteq\calV$ and $\calA'\subseteq\calV$, it is $f(\calA\cup \calA')\leq f(\calA)+f(\calA')$; ineq.~\eqref{aux1:4} holds for the same reason as ineq.~\eqref{aux1:2}; ineq.~\eqref{aux1:7} holds since or all elements $y \in \calY$ and all elements $p \in \calP$ it is $f(y)\geq f(p)$; finally, ineq.~\eqref{aux1:5} holds because $f$ is monotone and $y\in\calY$.
\hfill $\square$

To prove ineq.~\eqref{ineq:aux_3}, since it is  $\calB^\star_2\neq \emptyset$ (and, as a result, also $\calA^+_1\neq \emptyset$) and for all elements $a \in \calA^+_1$ and all elements $b\in \calB^\star_2$ it is $f(a)\geq f(b)$,  from Lemma~\ref{lem:D3} we have
\begin{align}
f(\calB^\star_2|\calA^+_1)&\leq |\calB^\star_2|f(\calA^+_1)\nonumber\\
&\leq \beta f(\calA^+_1),\label{aux:111}
\end{align}
since $|\calB^\star_2|\leq \beta$.  Overall,
\begin{align}
f(\calA^+_1)&\geq \frac{1}{\beta}f(\calB^\star_2|\calA^+_1)\label{aux5:1}\\
&\geq \frac{1}{\beta}f(\calB^\star_2|\calA^+_1\cup \calA^+_2)\label{aux5:2}\\
&=\frac{1}{\beta}f(\mathcal{B}_2^\star|\mathcal{A}_\textsc{Res}\setminus \mathcal{B}^\star(\mathcal{A}_\textsc{Res}))\label{aux5:3}\\
&=\eta\frac{1}{\beta}f(\calA_2),\label{aux5:4}
\end{align}
where ineq.~\eqref{aux5:1} to eq.~\eqref{aux5:4} hold for the following reasons: ineq.~\eqref{aux5:1} follows from ineq.~\eqref{aux:111}; ineq.~\eqref{aux5:2} holds since $f$ is submodular and $\calA_1^+\subseteq \calA_1^+\cup \calA_2^+$; eq.~\eqref{aux5:3} holds due to the definitions of the sets $\calA_1^+$, $\calA_2^+$ and $\mathcal{B}^\star(\mathcal{A}_\textsc{Res})$; finally, eq.~\eqref{aux5:4} holds due to the definition of $\eta$.  Overall, the~latter derivation concludes the proof of ineq.~\eqref{ineq:aux_3}.

\paragraph{Proof of ineq.~\eqref{ineq:aux_4}}  Let $b=1/\beta$.  We complete the proof first for the case where 
$(1-\eta)\geq \eta b$, and then for the case $(1-\eta)<\eta b$: i) When $(1-\eta)\geq \eta b$, $\max\{(1-\eta),\eta b\}= 1-\eta$ and $\eta \leq 1/(1+b)$.  Due to the latter, $1-\eta \geq b/(1+b)=1/(\beta+1)$ and, as a result,~\eqref{ineq:aux_4} holds. ii) When $(1-\eta)< \eta b$, $\max\{(1-\eta),\eta b\}= \eta b$ and $\eta > 1/(1+b)$. Due to the latter, $\eta b >  b/(1+b)$ and, as a result,~\eqref{ineq:aux_4} holds. 

We completed the proof of~$0\leq \eta\leq 1$, and ineqs.~\eqref{ineq:aux_2},~\eqref{ineq:aux_3} and~\eqref{ineq:aux_4}.  Thus, we also completed the proof of ineq.~\eqref{ineq:optimality_result_second_case_2}, that is, the last part of the proof of Algorithm~\ref{alg:rob_sub_max}'s approximation bound for $\kappa_f\neq 0$.

\subsection{Proof of Algorithm~\ref{alg:rob_sub_max}'s running time} We complete the proof in two steps, where the time for each evaluation of $f$ is denoted as $\tau_f$: we~first compute the time that line~2 of Algorithm~\ref{alg:rob_sub_max} needs to be executed, and then the time that lines 4-7 need to be executed. \mbox{Line 2} needs $m\log(m)+m\tau_f+m+O(\log(m))$ time, since it asks for $m$ evaluations of $f$, and their sorting, which takes $m\log(m)+m+O(\log(m))$ time, using, e.g., the merge sort algorithm.  Lines 4-7 need $m(\alpha-\beta)\tau_f+(\alpha-\beta)m$ time, since they ask for at most $m(\alpha-\beta)$ evaluations of $f$, since the while loop is repeated $\alpha-\beta$ times, and during each loop at most~$m$ evaluations of $f$ are needed by line~5, and for a maximal element in line~5, which needs $m$ time to be found.  Overall, Algorithm~\ref{alg:rob_sub_max} runs in $m(\alpha-\beta)\tau_f+(\alpha-\beta)m+m\log(m)+m\tau_f+m+O(\log(m))=O(m(\alpha-\beta)\tau_f)$ time.



\bibliographystyle{IEEEtran}
\bibliography{references,security_references}

\end{document}